\documentstyle[12pt]{article}

\def\fr{\frac}

\def\p{\partial}

\def\al{\alpha}
\def\Si{\Sigma}

\def\de{\delta}

\def\ha{\fr{1}{2}}

\def\a{\bar{a}}
\def\e{\exp_q}

\newcommand{\bd}{\begin{displaymath}}
\newcommand{\ed}{\end{displaymath}}
\newcommand{\bb}{\begin{equation}}
\newcommand{\ee}{\end{equation}}

\begin{document}
\baselineskip 1.3 \baselineskip

\vspace{.2cm}

\begin{center}
\large {\bf Some Realization of $gl_q(n)$-covariant Oscillator Algebra 
and $gl_q(n)$-covariant q-Virasoro Algebra with $q$ a root of unity.}
\\[1cm]
\large W.-S.Chung \\[.3cm]
\normalsize  
Department of Physics and Research Institute of Natural Science, \\
\normalsize College of Natural Sciences,  \\
\normalsize  Gyeongsang National University,   \\
\normalsize   Jinju, 660-701, Korea
\end{center}

\vspace{0.5cm}
\begin{abstract}
In this paper some realization of $gl_q(n)$-covariant oscillators is
obtained when $q$ is a root of unity. 
And the $gl_q(n)$-covariant q-Virasoro algebra is presented by using the
$gl_q(n)$-covariant
oscillators.
\end{abstract}

\section{Introduction}

Quantum                                                               
groups                                                                      
or                        
q-deformed                    
Lie         
algebra                
implies   
some    
specific 
deformations 
of                    classical                      Lie           algebras.

From           
a             
mathematical      
point   
of  
view,     
it 
is  
a 
non-commutative                                                  
associative                                                                 
Hopf         
algebra.         
The       
structure        
and 
representation 
theory    of                                   quantum               groups  
have                 
been             
developed        
extensively  
by  
Jimbo    
[1] 
and 
Drinfeld                                                                
[2].                                                                        
            
The                                   
q-deformation                    
of        
Heisenberg       
algebra   
was     
made  
by 
Arik and Coon [3], Macfarlane [4] and Biedenharn [5].
Recently                                                               
there                                               
has                           
been                    
some              
interest               
in    
more     
general   
deformations 
involving                                                                 
an                                                                          
arbitrary                      
real                           
functions     
of           
weight    
generators  
and   
including 
q-deformed algebras as a special case [6-10].

In this paper we review the $gl_q(n)$ oscillator algebra [11,12] and its
representation
when $q$ is a root of unity.
In this case we obtain some realizations of the algebra. We use these results
to present 'new' extension of q-Virasoro algebra which I call the
$gl_q(n)$-covariant
q-Virasoro algebra.

\section{$gl_q(n)$-Covariant Oscillator Algebra}

$gl_q(n)$-covariant oscillator algebra is defined as [11,12]
\bd
\a_i\a_j =\sqrt{q} \a_j \a_i,~~~(i<j)  
\ed
\bd
a_ia_j=\fr{1}{\sqrt{q}}a_j a_i,~~~(i<j)
\ed
\bd
a_i\a_j=\sqrt{q} \a_ja_i,~~~(i \neq j)
\ed
\bd
a_i\a_i =1+q \a_ia_i +(q-1) \Si_{k=i+1}^n\a_k a_k,~~~(i=1,2,\cdots,n-1)
\ed
\bd
a_n \a_n =1+q \a_n a_n,
\ed
\bb
[N_i, a_j]=-\de_{ij}a_j,~~~[N_i, \a_j]=\de_{ij}\a_j,~~~(i,j=1,2,\cdots, n )
\ee

Here $N_i$ plays a role of number operator and $a_i(\a_i)$ plays a role of
annihilation(creation) operator.
From now on, we restrict our concern to the case that $q$ is an N-th
primitive root
of unity, $ q=e^{2\pi i/N}$.
From the
above algebra one can obtain the relation between the number operators
and mode opeartors as follows 
\bb
\a_ia_i=q^{\Si_{k=i+1}^nN_k}[N_i],
\ee                                                      

where $[x]$ is called a q-number and is defined as
\bd
[x]=\fr{q^{x}-1}{q-1}.
\ed
\def\nn{|n_1,n_2,\cdots,n_n>}
Let us introduce the Fock space basis $\nn$ for the number operators
$N_1,N_2,\cdots, N_n$ satisfying
\bb
N_i\nn=n_i\nn,~~~(n_1,n_2,\cdots,n_n=0,1,2\cdots)
\ee
Then we have the following representation
\bd
a_i\nn=\sqrt{q^{\Si_{k=i+1}^nn_k}[n_i]}|n_1,\cdots, n_i-1,\cdots,n_n>
\ed
\bb
\a_i\nn=\sqrt{q^{\Si_{k=i+1}^nn_k}[n_i+1]}|n_1,\cdots, n_i+1,\cdots,n_n>.
\ee

From the above representation we know that there exists the ground state
$|0,0,\cdots,0>$ satisfying
$a_i|0,0>=0$ for all $i=1,2,\cdots,n$. Thus the state $\nn$ is obtatind by
applying the creation operators
to the ground state $|0,0,\cdots,0>$
\bb
\nn=\fr{\a_n^{n_n}\cdots\a_1^{n_1}}{\sqrt{[n_1]!\cdots [n_n]!}}|0,0,\cdots,0>.
\ee
If we introduce the scale operators as follows
\bb
Q_i=q^{N_i},~~(i=1,2,\cdots,n),
\ee
we have from the algebra (1)
\bb
[a_i,\a_i]=Q_iQ_{i+1}\cdots Q_n.
\ee
Acting the operators $Q_i$'s on the basis $\nn$ produces
\bb
Q_i\nn=q^{n_i}\nn .
\ee

In this representation the following relations are true:
\bb
[a_i,(\a_i)^n]=q^{1-n}[n]Q_i\cdots Q_n (\a_i)^{n-1}
\ee
We must notice that if $q$ is a N-th primitive root of unity, the
dimension of the representation space becomes finite.
In this case we obtain the following operator identities;
\bb
a_i^N=0~~~~(\a_i)^N=1
\ee

\section{First Realization}

The q-deformation of the Bargmann-Fock reptesentation of $gl_q(n)$
oscillators is realized by going over the space of analytic function
of n complex variables $z_1,\cdots,z_n$ such as $ |z_i|^2 \leq
(1-q)^{-1}$, where the operators $a_i,\a_i,N_i$ are defined
\bd
\a_i=T_{i+1}^{1/2}\cdots T_n^{1/2} z_i
\ed
\bd
a_i=\fr{1}{z_i}T_{i+1}^{1/2}\cdots T_n^{1/2} \fr{1-T_i}{1-q}
\ed
\bb
N_i=z_i=\fr{\p}{\p z_i}.
\ee
Here $T_i$ is a scaling operator satisfying
\bd
T_i^{\al}f(z_1,\cdots,z_i,\cdots,z_n)=f(z_i,\cdots,q^{\al}z_i,\cdots,z_n)
\ed
and the relation between $T_i$ and $N_i$ is
\bd
T_i=q^{N_i}
\ed
It can be easily checked that (11) satisfy the $gl_q(n)$ covariant
oscillator algebra (1). 

In the space of analytic n-variable function $f(z_1,\cdots,z_n)$ ,
there exists an inner product such as
\bb
(f,g)=\int (\Pi_{i=1}^n d_q^2 z_i)\mu(z_1,\cdots,z_n)
\bar{f}(z_1,\cdots,z_n) g(z_1,\cdots,z_n)
\ee
where the measure function $\mu(z_1,\cdots,z_n) $is defined as
\bb
\mu(z_1,\cdots,z_n)=\fr{1}{\pi^n}\Pi_{i=1}^n(exp_q(q|z_i|^2))^{-1}
\ee
where q-exponential function is defined as
\bd
\e(x)=\Si_{n=0}^{\infty}\fr{x^n}{[n]!}.
\ed
The q-exponential function satisfies the following recurrence relation
\bb
\e(q x)=[1-(1-q)x]\e(x)
\ee

From the above definition of the inner product, one can obtain the
orthonormal basis in the space of analytic n-variable function as
follows
\bb
u_n(\bar{z}_1,\cdots,\bar{z}_n)=\fr{\Pi_{i=1}^n \bar{z}_i^{n_i}}{\sqrt
{\Pi_{i=1}^n[n_i]!}}
\ee
Then we can easily see that
\bb
(u_n,u_m)=\de_{nm}
\ee
where we used
\bb
\int_0^{1/1-q}x^n e_q(qx)^{-1}d_qx =[n]!
\ee

\section{Second Realization}
Let $S_N$ be the discretized circle with points only at the positions
of N-th primitive root of unity,
\bd
S_N=\{1,q,q^2,\cdots,q^{N-1}\}.
\ed
Then the basis function of $gl_q(n)$ covariant oscillator system becomes a
$n$ variable function $f(q^{n_1},q^{n_2},\cdots,q^{n_n})$ defined on
the domain $S_N\otimes S_N\otimes \cdots \otimes S_N$ ( n-copies).
Now we define the action of the exponetials of the position and
momentum operators
\bd
(h_i f)(q^{n_1},\cdots,q^{n_i},\cdots,q^{n_n})
=f(q^{n_1},\cdots,q^{n_i-1},\cdots,q^{n_n}),
\ed
\bd
(g_i f)(q^{n_1},\cdots,q^{n_i},\cdots,q^{n_n})
= q^{n_i}f(q^{n_1},\cdots,q^{n_i},q^{n_{i+1}-1}\cdots,q^{n_n-1})
\ed
where $ n_i=0,1,\cdots,N-1$.
When $q$ is a N-th primitive root of unity , they satisfy
\bd
h_ig_i=qg_ih_i,
\ed
\bd
g_ig_j=qg_jg_i~~~~(i<j),
\ed
\bd
g_i^N=h_i^N=1
\ed
\bd
h_ig_j=g_jh_i~~~~(i \neq j)
\ed
\bb
h_ih_j=h_jh_i
\ee            
\def\a{\bar{a}}
Then $a_i$ and $\a_i$ operators in algebra (1) are realized as follows
\bd
\a_i =g_i,
\ed
\bb
a_i=g_i^{-1}\fr{h_{i+1}^2 \cdots h_n^2 (1-h_i)^2 }{1-q}
\ee
It can be easily checked that the above realization obeys the $gl_q(n)$
covariant oscillator algebra. From the above realization we can easily see
that the $\a_i$'s are
idempotent while $a_i$'s are nilpotent;
\bd
(\a_i)^N=1
\ed
\bb
a_i^N=0
\ee

\section{$gl_q(n)$-Covariant q-Virasoro Algebra}

Finally we would like to to mention that the realization of the
$gl_q(n)$ oscillator algebra may be used to construct the
$gl_q(n)$-covariant q-Virasoro algebra.

Like the q-Virasoro generator [13], we introduce the $l_n^{(i)}$ as
follows
\bb
l_m^{(i)}=(\a_i)^{m+1}a_i
\ee
which lead to the algebra
\bd
[l_r^{(i)},l_m^{(i)}]=q^{-m}[m-r]Q_iQ_{i+1}\cdots Q_n l_{r+m}^{(i)}
\ed
\bb
l_r^{(i)}l_m^{(j)}=q^{\ha m(r+2)}l_m^{(j)}l_r^{(i)},~~~(i<j)
\ee
If we redefine the operators through
\bb
L_m^{(i)}=Q_i^{-1}\cdots Q_n^{-1}l_m^{(i)}
\ee
we then have the following algebra
\bd
q^{m-r}L_r^{(i)}L_m^{(i)}-L_m^{(i)}L_r^{(i)}=[m-r]L_{r+m}^{(i)} \ed
\bb
L_r^{(i)}L_m^{(j)}=q^{\ha mr}L_m^{(j)}L_r^{(i)},~~~(i<j)
\ee
which goes to the n-copies of the ordinary Virasoro algebra in the $q
\rightarrow 1 $ limit.
However, when $q$ is a root of unity, it is something new although we do not
know
its meaning and physical application.

\section{Concluding Remark}
In this paper we studied the $gl_q(n)$-covariant multimode oscillator
algebra and some
of its realization when the deformation parameter $q$ is an N-th primitive
root of unity.
In this realization we found that the q-creation operator is idempotent
while the
q-annihilation operator is nilpotent, so the representation space becomes
finite 
dimensional.
We used these results to construct the $gl_q(n)$-covariant q-Virasoro
generator and to
find the new kind of deformed Virasoro algebra which I call the
$gl_q(n)$-covariant 
q-Virasoro algebra.

\section*{Acknowledgement}
This                   paper                was
supported         by  
the   KOSEF (961-0201-004-2)   
and   the   present   studies    were   supported   by   Basic  
Science 
Research Program, Ministry of Education, 1995 (BSRI-95-2413).

\vfill\eject

\end{document}